\documentclass{article}

\usepackage{amsmath}
\usepackage{amssymb}
\usepackage{amsbsy}
\usepackage{amscd}
\catcode`\@=11\relax
\newdimen\p@renwd
\font\tenex=cmex10 \setbox0=\hbox{\tenex B} \p@renwd=\wd0
\def\bbordermatrix#1{\begingroup \m@th
\setbox\z@\vbox{\def\\{\crcr\noalign{\kern2\p@\global\let\cr\endline}}%
    \ialign{$##$\hfil\kern2\p@\kern\p@renwd&\thinspace\hfil$##$\hfil
      &&\quad\hfil$##$\hfil\crcr
      \omit\strut\hfil\crcr\noalign{\kern-\baselineskip}%
      #1\crcr\omit\strut\cr}}%
  \setbox\tw@\vbox{\unvcopy\z@\global\setbox\@ne\lastbox}%
  \setbox\tw@\hbox{\unhbox\@ne\unskip\global\setbox\@ne\lastbox}%
  \setbox\tw@\hbox{$\kern\wd\@ne\kern-\p@renwd\left[\kern-\wd\@ne
    \global\setbox\@ne\vbox{\box\@ne\kern2\p@}%
    \vcenter{\kern-\ht\@ne\unvbox\z@\kern-\baselineskip}\,\right]$}%
  \null\;\vbox{\kern\ht\@ne\box\tw@}\endgroup}
\catcode`\@=12\relax

\newcommand{\cc}{{\mathbb C}}
\newcommand{\ZZ}{{\mathbb Z}}

\newcommand{\bfh}{{\boldsymbol{h}}}

\newcommand{\boldf}{{\boldsymbol{f}}}
\newcommand{\x}{{\boldsymbol{x}}}

\newcommand{\zero}{{\boldsymbol{0}}}
\newcommand{\bfl}{{\boldsymbol{\lambda}}}
\newcommand{\bfp}{{\boldsymbol{\partial}}}

\newcommand{\IN}{{\operatorname{in}}}
\newcommand{\supp}{{\operatorname{supp}}}
\newcommand{\corank}{{\operatorname{corank}\,}}

\def\D0{D_\zero}
\def\Span{\operatorname{Span}}
\newcommand{\ignore}[1]{}

\newtheorem{theorem}{Theorem}[section]
\newtheorem{lemma}[theorem]{Lemma}
\newtheorem{proposition}[theorem]{Proposition}
\newtheorem{definition}[theorem]{Definition}
\newtheorem{example}[theorem]{Example}
\newtheorem{algorithm}[theorem]{Algorithm}
\newtheorem{corollary}[theorem]{Corollary}
\newtheorem{remark}[theorem]{Remark}
\newcommand{\bC}{{\mathbb C}}
\newcommand{\bZ}{{\mathbb Z}}

\newcommand\qed{{\hspace*{\fill}$\Box$\vskip12pt plus 1pt}}
\newcommand\qex{\hfill $\diamond$}  

\begin{document}

\title{Higher-Order Deflation for Polynomial Systems
       with Isolated Singular Solutions}

\author{
Anton Leykin\thanks{
Department of Mathematics, Statistics, and Computer Science,
University of Illinois at Chicago, 851 South Morgan (M/C 249),
Chicago, IL 60607-7045, USA.
{\em Email:} leykin@math.uic.edu.
{\em URL:} http://www.math.uic.edu/{\~{}}leykin.}
\and
Jan Verschelde\thanks{
Department of Mathematics, Statistics, and Computer Science,
University of Illinois at Chicago, 851 South Morgan (M/C 249),
Chicago, IL 60607-7045, USA.
{\em Email:} jan@math.uic.edu or jan.verschelde@na-net.ornl.gov.
{\em URL:} http://www.math.uic.edu/{\~{}}jan.
This material is based upon work
supported by the National Science Foundation under Grant No.\
0105739 and Grant No.\ 0134611.}
\and
Ailing Zhao\thanks{
Department of Mathematics, Statistics, and Computer Science,
University of Illinois at Chicago, 851 South Morgan (M/C 249),
Chicago, IL 60607-7045, USA.
{\em Email:} azhao1@uic.edu.
{\em URL:} http://www.math.uic.edu/{\~{}}azhao1.}
}

\date{4 January 2007}

\maketitle
\begin{abstract}

Given an approximation to a multiple isolated solution of a polynomial
system of equations, we have provided a symbolic-numeric deflation
algorithm to restore the quadratic convergence of Newton's method.
Using first-order derivatives of the polynomials in the system,
our method creates an augmented system of equations which has the
multiple isolated solution of the original system as a regular root.

In this paper we consider two approaches to computing the
``multiplicity structure'' at a singular isolated solution.
An idea coming from one of them gives rise to our new
higher-order deflation method.
Using higher-order partial derivatives of the original polynomials,
the new algorithm reduces the multiplicity faster than our first method
for systems which require several first-order deflation steps.

We also present an algorithm to predict the order of the deflation.



\end{abstract}

\newpage

\section{Introduction}

This paper describes a numerical treatment of singular solutions
of polynomial systems. A trivial example to consider 
would be a single equation with a double root, $f(x) = x^2 = 0$, 
or a cluster of two very close roots, $f(x) = x^2 - \varepsilon^2  = 0$, 
where $0<\varepsilon \ll \mbox{machine precision}$. 
In both cases getting good approximate solutions with straightforward 
numerical approaches such as Newton's method is not easy.
Instead of attempting to solve the given equations we replace them with 
the system augmented by the equation's derivative,
$\bar \boldf (x)=(f(x),f'(x))=\zero$. 
Note that this completely symbolic procedure leads to a system 
with exact regular root in the first case, whereas in the second 
the system $\bar \boldf (x)=\zero$ would be inconsistent. 
However, a numerical solver applied to the latter would converge 
to a regular solution of a close-by system.   

In general setting, given an overdetermined system of equations in
many variables with a multiple isolated solution (a cluster of solutions)
our approach deflates the multiplicity of the solution (cluster)
by applying a certain numerical procedure. 
From the point of view of the numerical analysis it may be called 
a \emph{reconditioning} method: to recondition means to
reformulate a problem so its condition number improves.

Our deflation method was first presented at~\cite{VZ04},
and then described in greater detail in~\cite{LVZ06}.
In~\cite{LVZ07}, a directed acyclic graph of Jacobian matrices was
introduced for an efficient implementation.

On input we consider clusters of approximate zeroes of systems
$F(\x) = (f_1(\x)$, $f_2(\x)$, $\ldots,$ $f_N(\x)) = \zero$ of $N$ equations
in $n$ unknowns $\x \in \cc^n$.  We assume the cluster approximates
an isolated solution~$\x^*$ of $F(\x) = \zero$.  Therefore, $N \geq n$.
As~$\x^*$ is a singular solution, the Jacobian matrix of $F(\x)$,
denoted by~$A(\x)$, is singular at~$\x^*$.
In particular, we have $r = {\rm Rank}(A(\x^*)) < n$.

In case $r = n-1$, consider a nonzero vector~$\bfl$ in the kernel
of~$A(\x^*)$, which we denote by $\bfl \in \ker(A(\x^*))$,
then the equations
\begin{equation}
  g_i(\x)= \sum_{i=1}^{n} \lambda_j
 \frac{\partial f_i(\x)}{\partial x_j}, \quad i=1,2,\ldots,N,
\end{equation}
vanish at~$\x^*$, because $r = {\rm Rank}(A(\x^*)) < n$.
For $r < n-1$, our algorithm reduces to the corank-1 case,
replacing $A(\x)$ by $A(\x) B$, where $B$ is a random complex
$N$-by-$(r+1)$ matrix.  For the uniqueness $\bfl \in \ker(A(\x^*))$,
we add a linear scaling equation $\langle \bfh ,\bfl \rangle = 1$
(using a random complex $(r+1)$-vector $\bfh$).
and consider the augmented system
\begin{equation}
   G(\x,\bfl) =
   \left\{
     \begin{array}{lcl}
        F(\x)        & = & \zero \\
        A(\x) B \bfl & = & \zero \\
        \multicolumn{1}{r}{\langle \bfh ,\bfl \rangle} & = & 1.
     \end{array}
   \right.
\end{equation}
Let us denote by $\mu_F(\x^*)$ the multiplicity of $\x^*$ as 
a solution of the system $F(\x)=0$. In~\cite{LVZ06} we proved 
that there is a $\bfl^*$ 
such that $\mu_G(\x^*,\bfl^*) < \mu_F(\x^*)$.
Therefore, our deflation algorithm takes at most $m-1$ stages to determine
$\x^*$ as a regular root of an augmented polynomial system.

\noindent {\bf Related work.}  The literature on Newton's method
is vast.  As stated in~\cite{Deu04}, Lyapunov-Schmidt reduction
(see also~\cite[\S 6.2]{Gov00}, \cite{ABHJ99}, \cite{Kun96},
and~\cite{Lij04})) stands at the beginning of every
mathematical treatment of singularities.
We found the inspiration to develop a symbolic-numeric deflation
algorithm in~\cite{OWM83}.  The symbolic deflation procedure
of~\cite{Lec02} restores the quadratic convergence of Newton's method
with a complexity proportional to the square of the multiplicity
of the root.  Algorithms to compute the multiplicity are presented
in~\cite{BPS06}, \cite{DZ05}, and~\cite{ST98}.

\noindent {\bf Our Contributions.} We establish the link between
two different objects describing what we call the
\emph{multiplicity structure} of an isolated singular solution:
the dual space of differential functionals and the
initial ideal with respect to a local monomial order, both
associated to the ideal generated by the polynomials in the system
in the polynomial ring.

Next, following the latter method, we explain how to compute a
basis of the dual space, first, following the ideas of Dayton and
Zeng~\cite{DZ05}, then using the approach of Stetter and
Thallinger~\cite{ST98}.
We provide a formal symbolic algorithm for each approach,
respectively called the DZ and ST algorithms;
the ingredients of the algorithms do not go beyond linear algebra.
Moreover, we present an algorithm to determine the order of
the deflation.

The formalism developed for DZ and ST algorithms found a natural
continuation in \emph{higher-order deflation} method that
generalizes and extends the first-order deflation in \cite{LVZ06}.
For the systems that require more than one deflation step by our
first algorithm, the new deflation algorithm is capable of
completing the deflation in fewer steps. 

\noindent {\bf Acknowledgements.}  The material in this paper
was presented by the first two authors at the workshops on
computational algebraic geometry and real-number complexity,
organized respectively by Teresa Krick \& Andrei Gabrielov and
Peter Buergisser \& Gregorio Malajovich.  The authors thank the
organizers for the opportunities to present their work at these
FoCM 2005 workshops.

\section{Statement of the Main Theorem \& Algorithms}
 
The matrices $A^{(d)}(\x)$ we introduce below coincide for $d=1$
with the Jacobian matrix of a polynomial system.  
They are generalizations of the Jacobian matrix, built along the
same construction as the matrices used in the computation of the
multiplicity by Dayton and Zeng in~\cite{DZ05}.

\begin{definition} {\rm
The \emph{deflation matrix $A^{(d)}(\x)$ of a polynomial system
$F = (f_1$, $f_2$, $\ldots,f_N)$} of $N$ equations in $n$ 
unknowns~$\x = (x_1,x_2,\ldots,x_n)$ 
is a matrix with elements in $\bC[\x]$.  
The rows of $A^{(d)}(\x)$ are indexed by $\x^\alpha f_j$, 
where $|\alpha| < d$ and $j=1,2,\ldots,N$.
The columns are indexed by partial differential operators
$\bfp^\beta = \frac{\partial^{|\beta|}}{\partial x_1^{\beta_1}
   \cdots \partial x_n^{\beta_n}}$, 
where $\beta\neq\zero$ and $|\beta|\leq d$.
The element at row $\x^\alpha f_j$ and column $\partial^\beta$
of $A^{(d)}(\x)$ is
\begin{equation}
   \partial^\beta\cdot(\x^\alpha f_j) =
   \frac{\partial^{|\beta|}(\x^\alpha f_j)}{\partial \x^\beta}.
\end{equation}
$A^{(d)}(\x)$ has $N_r$ rows and $N_c$ columns,
$N_{r} = N \cdot \binom{n+d-1}{n}$ and $N_{c} = \binom{n+d}{n} - 1$.
}
\end{definition}

\begin{example}[Second-order deflation matrix]{\rm
Consider a system of 3 equations in 2 variables
$F=(f_1,f_2,f_3)=\zero$, where $f_1 = x_1^2$, $f_2 = x_1^2 - x_2^3$,
and $f_3 = x_2^4$. 
Then the second-order deflation matrix $A^{(2)}(\x)$ of $f$ is
\begin{equation}
\bbordermatrix{
 & \partial_{x_1} & \partial_{x_2} & \partial_{x_1}^2 
 & \partial_{x_1} \partial_{x_2} & \partial_{x_2}^2 \\
      f_1 & 2 x_1     &    0         & 2     & 0        &  0         \\
      f_2 & 2 x_1     & -3 x_2^2     & 2     & 0        & -6 x_2     \\
      f_3 & 0         & 4 x_2^3      & 0     & 0        & 12 x_2^2   \\
  x_1 f_1 & 3 x_1^2   & 0   & 6 x_1  & 0     & 0                     \\
  x_1 f_2 & 3 x_1^2   & -3 x_1 x_2^2 & 6 x_1 & 3 x_2^2  & -6 x_1 x_2 \\
  x_1 f_3 & x_2^4     & 4 x_1 x_2^3  & 0     & 0        & 0          \\
  x_2 f_1 & 2 x_1 x_2 & x_1^2        & 2 x_2 & 2 x_1    & 0          \\
  x_2 f_2 & 2 x_1 x_2 & -4 x_2^3     & 2 x_2 & 2 x_1    & -12 x_2^2  \\
  x_2 f_3 & 0         & 5 x_2^4      & 0     & 0        & 20 x_2^3
 }.
\end{equation}
Notice that $A^{(1)}(\x)$ (or the Jacobian matrix of $F$)
is contained in the first three rows and two columns of~$A^{(2)}(\x)$.
}
\end{example}

\begin{definition} {\rm
Let $\x^*$ be an isolated singular solution of the system $F(\x) = \zero$
and let $d$ be the order of the deflation.
Take a nonzero $N_{r}$-vector $(\lambda_\beta)_{\beta\neq\zero,\
|\beta|\leq d}$ in the kernel of $A^{(d)}(\x^*)$. 
It corresponds to what we call a \emph{deflation operator} -- a
linear differential operator with constant coefficients $\lambda_\beta$
\begin{equation}
   Q = \sum_{\beta\neq\zero,\ |\beta|\leq d}
   \lambda_\beta\partial^\beta \ \in \bC[\bfp].
\end{equation}
We use $Q$ to define $N_{r}$ new equations
\begin{equation}
  g_{j,\alpha}(\x) = Q \cdot (\x^\alpha f_j) = 0,
  \quad j=1,2,\ldots,N,\ |\alpha|<d.
\end{equation}
When we consider $\lambda_\beta$ as indeterminate,
we write $g_{j,\alpha}$ as~$g_{j,\alpha}(\x,\bfl)$.
In that case, we define $m$ additional linear equations,
for $m = \corank(A^{(d)}(\x^*))$:
\begin{equation}
  h_k(\bfl) = \sum_\beta b_{k, \beta} \lambda_\beta - 1 = 0,
  \quad k = 1,2,\ldots,m,
\end{equation}
where the coefficients $b_{k,\beta}$ are randomly chosen complex numbers.
}
\end{definition}

Now we are ready to state our main theorem.

\begin{theorem} \label{thm-deflation-step-general}
Let $\x^*\in\bC^n$ be an isolated solution of $F(\x)=\zero$.
Consider the following system in $\bC[\x,\bfl]$:
\begin{equation}\label{equ-d-deflation-system}
  G^{(d)}(\x,\bfl) =
  \left\{
    \begin{array}{rcll}
       f_j(\x) & = & 0, & j=1,2,\ldots,N; \\
       g_{j,\alpha}(\x,\bfl) & = & 0, & j=1,2,\ldots,N,\ |\alpha|<d; \\
       h_k(\bfl) & = & 0, & k=1,2,\ldots,m.
    \end{array}
  \right.
\end{equation}
For a generic choice of coefficients $b_{k,\beta}$, there exists a
unique $\bfl^*\in\bC^{N_{c}}$ such that the system
$G^{(d)}(\x,\bfl)$ has an isolated solution at $(\x^*,\bfl^*)$.
Moreover, the multiplicity of $(\x^*,\bfl^*)$ in $G(\x,\bfl)=\zero$
is strictly less than the multiplicity of $\x^*$ in $F(\x)=\zero$.
\end{theorem}

To determine the order $d$, we propose
Algorithm~\ref{alg_MinOrderForCorankDrop}.
This $d$ is then used 
in Algorithm~\ref{alg_SymbolicNumericDeflation}.

\begin{algorithm} \label{alg_MinOrderForCorankDrop}
{\rm d = {\bf MinOrderForCorankDrop}($F,\x_0$)

\begin{tabular}{ll}
Input: & $F$ is a finite set of polynomials; \\
       &  $\x^0 \approx \x^*$,
          $\x^*$ is an isolated multiple solution of $F(\x)=\zero$. \\
Output: & $d$ is the minimal number such that the augmented system $G^{(d)}$ \\
        & produced via a generic deflation operator $Q$ of order~$d$ has \\
        & corank of the Jacobian at $\x^*$ lower than $\corank A(\x^*)$. \\
\end{tabular}

\begin{tabular}{l}
take a generic vector 
$\gamma =(\gamma_1,\ldots,\gamma_n) \in \ker A(\x^0)$; \\
let $H(t) = F(\x^0+\gamma t) = F(x_1^0+\gamma_1 t,\ldots,x_n^0+\gamma_n t)$; \\
$d := \min\{ \ a \ |\  a \mbox{\rm ~belongs~to~the~support~of~} H(t) \ \} - 1$.
\end{tabular} }
\end{algorithm}

\begin{algorithm} \label{alg_SymbolicNumericDeflation}
{\rm $D^{(d)}F$ = {\bf Deflate}($F$,$d$,$\x_0$)

\begin{tabular}{ll}
Input: & $F$ is a finite set of polynomials in $\bC[\x]$; \\
       & $d$ is the order of deflation; \\
       & $\x^0 \approx \x^*$, 
         $\x^*$ is an isolated multiple solution of $F(\x)=\zero$. \\
Output: & $D^{(d)}F$ is a finite set of polynomials in $\bC[\x,\bfl]$ \\
        & such that there is $\bfl^*$ with 
          $\mu_{D^{(d)}F}(\x^*,\bfl^*) < \mu_{F}(\x^*)$.
\end{tabular}

\begin{tabular}{l}
determine the numerical corank $m$ of $A^{(d)}(\x^0)$; \\
return ${D^{(d)}F} := G^{(d)}(\x,\bfl)$ 
as in (\ref{equ-d-deflation-system}). 
\end{tabular}
}
\end{algorithm}
 
\section{Multiplicity Structure} \label{sec-multiplicity}

This section relates two different ways to obtain the multiplicity
of an isolated solution, constructing its {\em multiplicity
structure}.  Note that by a ``multiplicity structure'' -- a term
without a precise mathematical definition -- we mean any structure
which provides more local information about the singular solution
in addition to its multiplicity.  In this section we mention two
different approaches to describe this so-called multiplicity
structure.

\begin{example}[Running example 1] \label{exa-running} {\rm
Consider the system
\begin{equation} \label{equ-running}
  F(\x)
  = \left\{
       \begin{array}{c}
          x_2^3 = 0 \\
          x_1^2 x_2^2 = 0 \\
          x_1^4 + x_1^3x_2 = 0.
       \end{array}
    \right.
\end{equation}
The system $F(\x) = \zero$ has only one isolated solution
at~$(0,0)$ of high multiplicity.  Below we will show how
to compute the multiplicity of~$(0,0)$.~\qex  }
\end{example}

\subsection{Standard Bases}

Assume $\zero\in\bC^n$ is an isolated solution
of the system $F(\x)=\zero$.
Let $I = \langle F \rangle \subset R =\bC[\x]$ be
the ideal generated by the polynomials in the system.
Given a local monomial order $\geq$, the initial ideal
$\IN_\geq(I) = \{ \IN_\geq(f)\ |\ f \in I \} \subset R$
describes the multiplicity structure of $\zero$ by means
of \emph{standard monomials}, i.e.:
monomials that are not contained in $\IN_\geq(I)$.
A graphical representation of a monomial ideal
is a monomial \emph{staircase}.

\begin{figure}[hbt]

\unitlength 0.8mm

\begin{picture}(200,50)(0,0)

\put(4,0){
\begin{picture}(80,50)(0,0)

\linethickness{0.3mm}
 \put(0,0){\vector(0,1){45}}
 \put(0,0){\vector(1,0){65}}
 \put(0,0){\circle{4}}
 \put(0,10){\circle{4}}
 \put(20,0){\circle{4}}
 \put(0,20){\circle{4}}
 \put(10,10){\circle{4}}
 \put(10,0){\circle{4}}
 \put(0,30){\circle*{4}}
 \put(40,0){\circle*{4}}
 \put(0,30){\line(1,0){20}}
\linethickness{0.6mm}
 \put(40,0){\line(0,1){20}}

\put(25,25){\makebox(0,0)[cl]{$\underline{x_1^2x_2^2}$}}
\put(44,4){\makebox(0,0)[cl]{$\underline{x_1^4}+x_1^3x_2$}}

\linethickness{0.6mm} \put(20,20){\line(0,1){10}}
\linethickness{0.6mm} \put(20,20){\line(1,0){20}}
\linethickness{0.3mm} \put(20,20){\circle*{4}}

\put(5,35){\makebox(0,0)[cl]{$\underline{x_2^3}$}}
\put(40,40){\makebox(0,0)[cl]{$w = (-1,-2)$}}

\linethickness{0.3mm}
\multiput(30,30)(0.12,0.24){83}{\line(0,1){0.24}}
\put(30,30){\vector(-1,-2){0.12}}

\linethickness{0.3mm} \put(30,0){\circle{4}}
\linethickness{0.3mm} \put(10,20){\circle{4}}
\linethickness{0.3mm} \put(20,10){\circle{4}}
\linethickness{0.3mm} \put(30,10){\circle{4}}

\end{picture}
}

\put(85,0){
\begin{picture}(80,50)(0,0)

\linethickness{0.3mm} 
 \put(0,0){\vector(0,1){45}}
 \put(0,0){\vector(1,0){60}}
 \put(0,0){\circle{4}}
 \put(0,10){\circle{4}}
 \put(20,0){\circle{4}}
 \put(0,20){\circle{4}}
 \put(10,10){\circle{4}}
 \put(10,0){\circle{4}}
 \put(0,30){\circle*{4}}
 \put(50,0){\circle*{4}}
 \put(0,30){\line(1,0){20}}
\linethickness{0.6mm}
\put(50,0){\line(0,1){10}}
\put(25,25){\makebox(0,0)[cl]{$\underline{x_1^2x_2^2}$}}

\put(66,14){\makebox(0,0)[cc]{}}
\put(35,15){\makebox(0,0)[cl]{$\underline{x_1^3x_2}+x_1^4$}}

\linethickness{0.6mm} \put(20,20){\line(0,1){10}}
\linethickness{0.6mm} \put(20,20){\line(1,0){10}}
\linethickness{0.3mm} \put(20,20){\circle*{4}}

\put(5,35){\makebox(0,0)[cl]{$\underline{x_2^3}$}}

\put(8,44){\makebox(0,0)[cl]{}}

\put(20,40){\makebox(0,0)[cl]{$w = (-2,-1)$}}

\put(0,20){\makebox(0,0)[cl]{}}


\linethickness{0.3mm}
\multiput(30,30)(0.24,0.12){83}{\line(1,0){0.24}}
\put(30,30){\vector(-2,-1){0.12}}
\linethickness{0.3mm} \put(30,0){\circle{4}}

\linethickness{0.3mm} \put(10,20){\circle{4}}
\linethickness{0.3mm} \put(20,10){\circle{4}}
\linethickness{0.3mm} \put(40,0){\circle{4}}
\linethickness{0.6mm} \put(30,10){\line(0,1){10}}
\linethickness{0.3mm} \put(30,10){\circle*{4}}
\linethickness{0.6mm} \put(30,10){\line(1,0){20}}
\put(55,5){\makebox(0,0)[cl]{$\underline{x_1^5}$}}

\end{picture}
}
\end{picture}

\caption{Two monomial staircases for two different monomial orderings
applied to the same system.
The full circles represent the generators of the initial ideals.
The multiplicity is the number of standard monomials,
represented by the empty circles under the staircase.}
\label{figstair}
\end{figure}
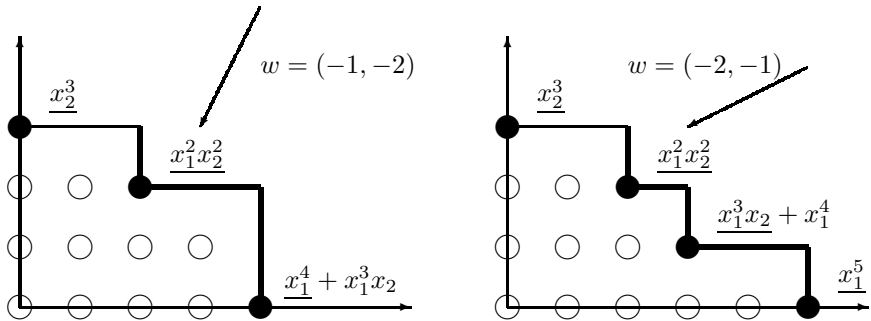

\begin{example}[Initial ideals with respect to a local order]
{\rm Consider the system~(\ref{equ-running})
of Example~\ref{exa-running}.

Figure~\ref{figstair} shows the staircases for initial ideals of
$I = \langle F \rangle$ w.r.t. two local weight orders $\geq_\omega$.
Computer algebra systems Macaulay~2~\cite{Macaulay2} 
and Singular~\cite{GPS01} can be used for these kind of computations, 
see also~\cite{GP96,GP02} for theory, in particular, on Mora's 
tangent cone algorithm~\cite{Mor82}.

In the example the leading monomials at the corners of the
staircase come from the elements of the corresponding \emph{standard basis}.
For the weight vector $w=(-1,-2)$ the original generators give such a 
basis (initial terms underlined).
For $w=(-2,-1)$ one more polynomial is needed.~\qex }
\end{example}

\subsection{Dual Space of Differential Functionals}

Another approach at the multiplicity structure is
described in detail in~\cite{Ste04,Tha98}; see also~\cite{Mou97}.
Using duality to define
the multiplicity of a solution goes back to Macaulay~\cite{Mac16}.
In this approach,
{\em differential functionals} are denoted by
\begin{equation}
   \Delta_\alpha(f) = \left.\frac{1}{\alpha_1!\cdots\alpha_n!}\cdot
   \frac{\partial^{|\alpha|}f}{\partial x^{\alpha_1}\cdots\partial
   x^{\alpha_n}}\right|_{\x=\zero}.
\end{equation}
Observe that
\begin{equation}
\ \Delta_\alpha(\x^\beta)
     = \left\{
          \begin{array}{ll}
            1, & \alpha=\beta \\
            0, & \alpha\neq\beta.
          \end{array}
       \right.
\end{equation}
We then define the local \emph{dual space}
of differential functionals $D_\zero[I]$ as
\begin{equation}
  D_\zero[I] = \{ L\in \Span\{\Delta_\alpha\ |\
                 \alpha\in\bZ_{\geq 0}^n\} |\
          L(f)=0 \mbox{ for all } f\in I\},
\end{equation}

\begin{example}[Dual space of running example 1] {\rm
For the ideal defined by the polynomials in
the system~(\ref{equ-running}) we have
\begin{equation}
\begin{array}{ccc}
  D_\zero[I] = \Span\{ &
   \underline{\Delta_{(4,0)}}-\Delta_{(3,1)},\ \underline{\Delta_{(3,0)}},\ \underline{\Delta_{(2,1)}},\ \underline{\Delta_{(1,2)}},& \\ 
  &\underline{\Delta_{(2,0)}},\ \underline{\Delta_{(1,1)}},\ \underline{\Delta_{(0,2)}},\ 
      \underline{\Delta_{(1,0)}},\ \underline{\Delta_{(0,1)}},\
    \underline{\Delta_{(0,0)}} &\}.
\end{array}
\end{equation}
Notice that here the basis of the dual space is chosen in such a
way that the (underlined) leading terms with respect to the weight
order $\geq_{(2,1)}$ correspond to the monomials under the
staircase in Example \ref{exa-running} for the order
$\geq_{(-2,-1)}$. We will show that it is not a coincidence later
in this section.\qex
}
\end{example}

\subsection{Dual Bases versus Standard Bases}

Since both local dual bases and initial ideals w.r.t. local orders
describe the same, there exists a natural correspondence between the two.

Let $\geq$ be an order on the nonnegative integer lattice
$\ZZ_{\geq 0}^n$ that defines a local monomial order and let
$\succeq$ be the opposite of $\geq$: i.e. $\alpha \succeq \beta
\Leftrightarrow \alpha \leq \beta$. (Note: $\succeq$ defines a
global monomial order.)

For a linear differential functional $L = \sum c_\alpha
\Delta_\alpha$ define the \emph{support}: $\supp(L) = \{
\alpha\in\ZZ_{\geq 0}^n \ |\ c_\alpha \neq 0 \}.$
For the dual space, $\supp (\D0[I]) = \bigcup_{L\in \D0[I]} \supp(L)$.

Using the order $\succeq$ we can talk about the leading or
\emph{initial term} of $L$: let $\IN_\succeq(L)$ be the maximal
element of $\supp(L)$ with respect to $\succeq$. Define the
\emph{initial support} of the dual space as $\IN_\succeq(\D0[I]) =
\{ \IN_\succeq(L)\ |\ L\in \D0[I]\}$. The initial support is
obviously contained in the support, in our running example the
containment is proper:
\begin{eqnarray*}
\IN_{(2,1)}(\D0[I])
   & = &  \{ (i,j)\ |\ i+j\leq 3 \} \cup \{ (4,0) \} \\
   & \subset & \{ (i,j)\ |\ i+j\leq 3 \} \cup  \{ (4,0) \cup (3,1) \} = \supp(\D0[I]).
\end{eqnarray*}

\begin{theorem}
The number of elements in the initial support equals the dimension
of the dual space, therefore, is the multiplicity. Moreover, with
the above assumptions on the orders $\geq$ and $\succeq$, the
standard monomials w.r.t. the local order $\geq$ are $\{\x^\alpha\
|\ \alpha \in \IN_\succeq(\D0[I]) \}$.
\end{theorem}

\noindent {\em Proof.}
Pick $L_\beta\in \D0[I], \beta\in\IN_\succeq(\D0[I])$ such that
$\IN_\succeq(L_\beta) = \beta$. One can easily show that
$\{L_\beta\}$ is a basis of $\D0[I]$.

Take a monomial $\x^\alpha \in \IN_\geq(I)$, then there is $f \in
I$ such that $\x^\alpha = \IN_\geq(f)$. Next, take any linear
differential functional $L$ with $\IN_\succeq(L) = \alpha$. Since
the orders $\geq$ and $\succeq$ are opposite, there are no similar
terms in the tail of $L$ and the tail of $f$, therefore, $L(f) =
\IN_\succeq(L)(\IN_\geq(f)) \neq 0$.

It follows that, $L\notin \D0[I]$, which proves that the set of
standard monomials is contained in the initial support of
$\D0[I]$. They are equal since they both determine the dimension.~\qed

Consider the ring of linear differential operators $\mathcal D=
\bC[\bfp]$ with the natural action (denoted by ``$\cdot$'') on
polynomial ring $R=\bC[\x]$.

\begin{lemma} \label{lem-drop}
Let $Q \in \bC[\bfp]$ and $f \in \bC[\x]$ such that
$\IN_\succeq(Q) \succeq \IN_\geq(f)$ (in $\ZZ_{\geq 0}^n$).

Then $\IN_\geq (Q\cdot f) = \IN_\geq(f) - \IN_\succeq(Q) \in
\ZZ_{\geq 0}^n$.
\end{lemma}

\section{Computing the Multiplicity Structure} \label{sec-computing-mult}

Let the ideal $I$ be generated by $f_1,f_2,\ldots,f_N$.
Let $\D0^{(d)}[I]$ the part of $\D0[I]$ containing functionals of
order at most $d$. We would like to have a criterion that for the
differential functional $L$ of degree at most $d$ guarantees $L\in
\D0^{(d)}[I]$.

Below we describe two such criteria referred to as
\emph{closedness conditions}; their names are arranged to match
the corresponding computational techniques of
Dayton-Zeng~\cite{DZ05} and Stetter-Thallinger~\cite{ST98}
that we will describe later respectively as the DZ and ST algorithms.

A functional $L = \sum c_\alpha \Delta_\alpha$ with $c_\alpha \in
\bC$ of order $d$ belongs to the dual space $\D0[I]$ if and only if

\begin{itemize}

\item \textbf{(DZ-closedness)} $L(g\cdot f_i) = 0$
      for all $i = 1,2,\ldots,N$ and polynomials $g(\x)$
      of degree at most $d-1$.

\item \textbf{(ST-closedness)} $L(f_i) = 0$ for all $i$ and
      $\sigma_j(L) \in \D0[I]$ for all $j = 1,2,\ldots,n$,
      where $\sigma_j : \D0[I] \to \D0[I]$ is a linear map such that
\begin{equation}
  \sigma_j(\Delta_\alpha)
   = \left\{ \begin{array}{cl}
                0, &\mbox{ if } \alpha_j=0,\\
                \Delta_{\alpha-e_j}, &\mbox{ otherwise}.
             \end{array}
     \right.
\end{equation}

\end{itemize}

The basic idea of both DZ and ST algorithms is the same:
build up a basis of $\D0$ incrementally by computing $\D0^{(d)}$ for
$d=1,2,\ldots$ using the corresponding closedness condition.
The computation stops when $\D0^{(d)} = \D0^{(d-1)}$.

\begin{example}[Running example 2]\label{exa-small}
{\rm
Consider the system in $\bC[x_1,x_2]$ given by three polynomials
$f_1 = x_1x_2$, $f_2 = x_1^2-x_2^2$, and $f_3 = x_2^4$,
which has only one isolated root at~$(0,0)$.~\qex
}
\end{example}

\subsection{The Dayton-Zeng Algorithm}

We shall outline only a summary of this approach,
see~\cite{DZ05} for details.

If $\zero$ is a solution of the system,
then $\D0^{(0)} = \Span\{\Delta_\zero\}$.

At step $d>0$, we compute $\D0^{(d)}$.
Let the functional
\begin{equation}
   L = \sum_{|\alpha|\leq d,\ \alpha\neq \zero} c_\alpha
   \Delta_\alpha
\end{equation}
belong to the dual space $\D0^{(d)}$.
Then the vector of coefficients $c_\alpha$ is in the kernel of the
following matrix $M_{DZ}^{(d)}$ with $N B(d-1)$ rows and $B(d)-1$
columns, where $B(d)=\binom{n+d}{n}$ is the number of monomials in
$n$ variables of degree at most $d$.

The rows of $M_{DZ}^{(d)}$ are labelled with $x^\alpha f_j$, where
$|\alpha| < d$ and $j=1,2,\ldots,N$.
The columns correspond to
$\Delta_\beta$, where $\beta\neq\zero,\ |\beta|\leq d$.
\begin{equation}
\mbox{[The entry of $M_{DZ}^{(d)}$ in row $x^\alpha f_j$ and
column $\Delta_\beta$]} = \Delta_\beta(x^\alpha f_j).
\end{equation}

At the step $d=3$ we have the following $M_{DZ}^{(3)}$
$${\small
\begin{array}{c||cc|ccc|cccc}
 & \Delta_{(1,0)} & \Delta_{(0,1)} & \Delta_{(2,0)}
 & \Delta_{(1,1)} & \Delta_{(0,2)} & \Delta_{(3,0)}
 & \Delta_{(2,1)} & \Delta_{(1,2)} & \Delta_{(0,3)} \\ \hline \hline
f_1 & 0 & 0 & 0 & 1 &  0 & 0 & 0 & 0 & 0  \\
f_2 & 0 & 0 & 1 & 0 & -1 & 0 & 0 & 0 & 0  \\
f_3 & 0 & 0 & 0 & 0 &  0 & 0 & 0 & 0 & 0  \\ \hline
x_1f_1 & 0 & 0 & 0 & 0 & 0 & 0 & 1 & 0 & 0 \\
x_1f_2 & 0 & 0 & 0 & 0 & 0 & 1 & 0 & -1 & 0 \\
x_1f_3 & 0 & 0 & 0 & 0 & 0 & 0 & 0 & 0 & 0 \\
x_2f_1 & 0 & 0 & 0 & 0 & 0 & 0 & 0 & 1 & 0 \\
x_2f_2 & 0 & 0 & 0 & 0 & 0 & 0 & 1 & 0 & -1 \\
x_2f_3 & 0 & 0 & 0 & 0 & 0 & 0 & 0 & 0 & 0 \\ \hline
x_1^2f_1  & 0 & 0 & 0 & 0 & 0 & 0 & 0 & 0 & 0   \\
\vdots & \vdots & \vdots & \vdots & \vdots & \vdots
& \vdots & \vdots & \vdots  & \vdots
\end{array}
}
$$
Note that the last block of 9 rows is entirely zero.

Analyzing the kernel of this matrix one sees that there are no
functionals of degree 3 in the dual space, which is then is equal
to $\D0^{(2)}[I]$
\begin{equation}
\D0[I] = \Span\{\Delta_{(0,0)}, \Delta_{(1,0)}, \Delta_{(0,1)},
\Delta_{(2,0)} + \Delta_{(0,2)}\}.
\end{equation}

\subsection{The Stetter-Thallinger Algorithm}

The matrix $M_{ST}^{(d)}$ is a matrix consisting of $n+1$ blocks
stacked on top of each other:

\begin{itemize}
\item The top block contains the first $N$ rows of $M_{DZ}^{(d)}$;

\item For every $j=1,2,\ldots,n$, let $S_j^{(d)}$ be the
$(B(d-1)-1)\times (B(d)-1)$-matrix for the linear map
\begin{equation}
  \sigma_j : \D0^{(d)}/\Span\{\Delta_\zero\} \to
\D0^{(d-1)}/\Span\{\Delta_\zero\}
\end{equation}
w.r.t. standard bases of functionals.

The block $M_{ST}^{(d-1)}S_j$ represents the closedness condition
for the ``anti-derivation'' $\sigma_j$.
\end{itemize}

Let us go through the steps of the algorithm for the Example
\ref{exa-small}.

\noindent \textbf{Step 1.} At the beginning we have $M_{ST}^{(1)}$
equal to
$$
\begin{array}{c||cc}

    & \Delta_{(1,0)} & \Delta_{(0,1)}\\
\hline \hline
f_1 & 0              & 0             \\
f_2 & 0              & 0             \\
f_3 & 0              & 0             \\

\end{array}
$$
Therefore, $\D0^{(1)} = \Span \{ \Delta_{(0,0)}, \Delta_{(1,0)},
\Delta_{(0,1)} \}$.

\noindent \textbf{Step 2.} Since $M_{ST}^{(1)}S_j^{(2)}=0$ for all
$j$, the matrix $M_{ST}^{(2)}$ is
$$
\begin{array}{c||cc|ccc}

    & \Delta_{(1,0)} & \Delta_{(0,1)} & \Delta_{(2,0)} & \Delta_{(1,1)} & \Delta_{(0,2)}\\
\hline \hline
f_1 & 0 & 0 & 0 & 1 & 0  \\
f_2 & 0 & 0 & 1 & 0 & -1 \\
f_3 & 0 & 0 & 0 & 0 & 0  \\ \hline
    & 0 & 0 & 0 & 0 & 0  \\
\vdots & \vdots & \vdots & \vdots & \vdots & \vdots
\end{array}
$$
Therefore, $\D0^{(2)} = \Span \{ \Delta_{(0,0)}, \Delta_{(1,0)},
\Delta_{(0,1)}, \Delta_{(2,0)} + \Delta_{(0,2)} \}$.

We can ``prune'' the matrix $M_{ST}^{(2)}$ by row-reducing it to
the following matrix with the same kernel:
$$
\tilde M_{ST}^{(2)} = \left[
\begin{array}{ccccc}
0 & 0 & 0 & 1 & 0  \\
0 & 0 & 1 & 0 & -1 \\
\end{array}
\right]
$$

\noindent \textbf{Step 3.} Compute $S_1^{(3)}$ that represents
$\sigma_1$:
$${\small
\begin{array}{c||cc|ccc|cccc}

               & \Delta_{(1,0)} & \Delta_{(0,1)} & \Delta_{(2,0)} & \Delta_{(1,1)} & \Delta_{(0,2)} & \Delta_{(3,0)} & \Delta_{(2,1)} & \Delta_{(1,2)} & \Delta_{(0,3)}\\
\hline \hline
\Delta_{(1,0)} & 0              & 0              & 1              & 0              & 0              & 0              & 0              & 0              & 0               \\
\Delta_{(0,1)} & 0              & 0              & 0              & 1              & 0              & 0              & 0              & 0              & 0               \\
\hline
\Delta_{(2,0)} & 0              & 0              & 0              & 0              & 0              & 1              & 0              & 0              & 0               \\
\Delta_{(1,1)} & 0              & 0              & 0              & 0              & 0              & 0              & 1              & 0              & 0               \\
\Delta_{(0,2)} & 0              & 0              & 0              & 0              & 0              & 0              & 0              & 1              & 0               \\
\end{array}
}
$$
The matrix $S_2^{(3)}$ can be defined similarly.

The top block of the matrix $M_{ST}^{(3)}$ is
$${\small
\begin{array}{c||cc|ccc|cccc}

    & \Delta_{(1,0)} & \Delta_{(0,1)} & \Delta_{(2,0)} & \Delta_{(1,1)} & \Delta_{(0,2)} & \Delta_{(3,0)} & \Delta_{(2,1)} & \Delta_{(1,2)} & \Delta_{(0,3)}\\
\hline \hline
f_1 & 0              & 0              & 0              & 1              & 0                         & 0              & 0              & 0              & 0               \\
f_2 & 0              & 0              & 1              & 0              & -1                        & 0              & 0              & 0              & 0               \\
f_3 & 0              & 0              & 0              & 0              & 0                         & 0              & 0              & 0              & 0               \\
\end{array}
}
$$

Despite the last 4 columns being 0, there are no new elements of
order 3 in the dual space due to the other two blocks: $ \tilde
M_{ST}^{(2)} S_1^{(3)}$:
$${\small
\begin{array}{c||cc|ccc|cccc}
   & \Delta_{(1,0)} & \Delta_{(0,1)} & \Delta_{(2,0)} & \Delta_{(1,1)}
   & \Delta_{(0,2)} & \Delta_{(3,0)} & \Delta_{(2,1)} & \Delta_{(1,2)}
   & \Delta_{(0,3)}\\
\hline \hline
x_1f_1 & 0              & 0              & 0              & 0              & 0              & 0              & 1              & 0              & 0               \\
x_1f_2 & 0              & 0              & 0              & 0              & 0              & 1              & 0              & -1             & 0               \\
\end{array}
}
$$
and $ \tilde M_{ST}^{(2)} S_2^{(3)}$:
$${\small
\begin{array}{c||cc|ccc|cccc}
       & \Delta_{(1,0)} & \Delta_{(0,1)} & \Delta_{(2,0)} & \Delta_{(1,1)} & \Delta_{(0,2)} & \Delta_{(3,0)} & \Delta_{(2,1)} & \Delta_{(1,2)} & \Delta_{(0,3)}\\
\hline \hline
x_2f_1 & 0              & 0              & 0              & 0              & 0              & 0              & 0              & 1              & 0               \\
x_2f_2 & 0              & 0              & 0              & 0              & 0              & 0              & 1              & 0              & -1               \\
\end{array}
}
$$

\medskip

Comparing to DZ algorithm, in step 3, we managed to avoid the
computation of 9 last zero rows of $M_{DZ}^{(3)}$ in this
particular example.  We now also see how its 4 last nonzero rows
show up in the ``closedness condition'' blocks of $M_{DZ}^{(3)}$.

\section{Proofs and Algorithmic Details}\label{sec-general-deflation}

In this section we justify the main theorems stated before and
give details about the algorithms presented above.

\subsection{First-Order Deflation}\label{sec-old-deflation}

In this section we summarize our deflation method introduced in
\cite{LVZ06}. Not only it is done for the convenience of the reader,
but also for our own convenience as we plan to build a
higher-order deflation algorithm in Section
\ref{sec-general-deflation} using the algorithm following the
pattern established in this section.

\medskip

\noindent \textbf{One deflation step with fixed $\bfl$.} The basic
idea of the method is relatively simple. Let $\bfl\in\bC^n$ be a
nonzero vector in $\ker(A(\x^*))$, then the equations
\begin{equation}
   g_i(\x)= \bfl\cdot\nabla f_i(\x) = \sum_{i=1}^{n} \lambda_j
\frac{\partial f_i(\x)}{\partial x_j}, \quad i=1,2,\ldots,N
\end{equation}
have $\x^*$ as a solution. Moreover,

\begin{theorem}\label{thm-fixed-lambda}
The augmented system
\begin{equation}
G(\x) = (f_1,\ldots,f_N,g_1,\ldots,g_N)(\x)=\zero
\end{equation}
of equations in $\bC[\x]$ is a {{\em deflation}} of the original
system $F(\x)=\zero$ at $x^*$, i.e.  $G(\x^*)=0$ and the
multiplicity of the solution $x^*$ is lower in the new system.
\end{theorem}

The original proof of this statement in \cite{LVZ06} uses the notion
of a standard basis of the ideal $I = (f_1,f_2,\ldots,f_N)$ in the
polynomial ring $R = \bC[\x]$ w.r.t. a local order; this tool of
computational commutative algebra can be used to obtain the
multiplicity of $\x^*$, which is defined as the $\bC$-dimension of
the local quotient ring $R_\x / R_x I$.

On the other hand it is in correspondence with another way of
looking at multiplicities -- dual spaces of local functionals, so
the proof can be written in that language as well
(see Section \ref{sec-multiplicity}).

\noindent \textbf{One deflation step with indeterminate $\bfl$.}
Without loss of generality, we may assume $\corank(A(\x^*))=1$;
consult \cite{LVZ06} to see how the general case is reduced to this.
Consider $N+1$ additional polynomials in $\bC[\x,\bfl]$ in $2n$
variables:
\begin{equation}
g_i(\x,\bfl) = \bfl\cdot\nabla f_i(\x) = \sum_{j=1}^{n} \lambda_j
\frac{\partial f_i(\x)}{\partial x_j},\ (i=1,2,\ldots,N)
\end{equation}
\begin{equation}
h(\bfl) = \sum_{j=1}^{n} b_j \lambda_j - 1,
\end{equation}
where the coefficients $b_j$ are random complex numbers.

\begin{theorem} \label{thm-deflation-step}
Let $\x^*\in\bC^n$ be an isolated solution of $F(\x)=\zero$ (in
$\bC[\x]$).

For a generic choice of coefficients $b_j$, $j=1,2,\ldots,n$, there
exists a unique $\bfl^*\in\bC^n$ such that the system
\begin{equation}
G(\x,\bfl) = (f_1,\ldots,f_N,g_1,\ldots,g_N,h)(\x,\bfl)=0
\end{equation}
of equations in $\bC[\x,\bfl]$ has an isolated solution at
$(\x^*,\bfl^*)$.

The multiplicity of $(\x^*,\bfl^*)$ in $G(\x,\bfl)=\zero$ is lower
than that of $\x^*$ in $F(\x)=\zero$.
\end{theorem}

\noindent {\em Proof.} Follows from Proposition 3.4 in \cite{LVZ06}.~\qed

Theorem \ref{thm-deflation-step} provides a recipe for the
deflation algorithm: one simply needs to keep deflating until the
solution of the augmented system corresponding to $\x^*$ becomes
regular.

As a corollary we have that the number of deflations needed to
make a singular isolated solution $\x^*$ regular is less than the
multiplicity of~$\x^*$.

\subsection{Higher-Order Deflation with Fixed Multipliers}

We use the deflation operator to define an augmented system.

\begin{theorem}\label{thm-fixed-lambda-general}
Let $f_1,f_2,\ldots,f_N$ form a standard basis of $I$ w.r.t. the order
opposite to $\succeq$.  Consider the system $G^{(d)}(\x) = \zero$
in $\bC[\x]$, where
\begin{equation}
G^{(d)}(\x) = \left\{\begin{array}{rc}
f_j(\x) &(j=1,2,\ldots,N)\\
g_{j,\alpha}(\x) &(j=1,2,\ldots,N,\ |\alpha|<d)\\
\end{array}\right..
\end{equation}
\begin{description}
\item{(a)} The system $G^{(d)}(\x) = \zero$ is a deflation of the
original system $F(\x)=\zero$ at~$\x^*$.
\item{(b)} Let $I = (F)$ and $J = (G^{(d)})$ be the ideals generated by
polynomials of the systems and $\succeq$ be a global monomial
order on $\ZZ_{\geq 0}^n$.
Then the following relation holds for initial supports
\begin{equation}
\IN_\succeq(\D0[J]) \subset \{\beta-\beta_Q\ |\
\beta\in\IN_\succeq(\D0[I])\} \cap \ZZ_{\geq 0}^n,
\end{equation}
where $\beta_Q$ is the maximal element of the set
$\IN_\succeq(\D0[I])\cap\{\beta : |\beta|\leq d \}$.
\end{description}

\end{theorem}

\noindent {\em Proof.}
Let $\bfl \in \ker (A^{(d)}(\x^*))$ be the vector used above to
construct the operator $Q \in \bC[\bfp]$ and the equations
$g_{j,\alpha}(\x)=0$.

First of all, $g_{j,\alpha}(\x^*) = (Q\cdot (\x^\alpha
f_j))|_{\x=\x^*} = 0$ provided $|\alpha| < d$ (by construction),
hence, $\x^*$ is a solution to $G^{(d)}(\x) = \zero$.

To prove (a), it remains to show that the multiplicity drops,
which follows from part (b) that is treated in the rest of this proof.

We shall assume for simplicity that $\x^*=\zero$. This is done
without the loss of generality using a linear change of
coordinates: $\x\mapsto\x+\x^*$. It is important to note that in
the new coordinates polynomials $Q\cdot(\x^\alpha f_j(\x+\x^*))$
generate the same ideal as the polynomials
$Q\cdot((\x-\x^*)^\alpha f_j(\x+\x^*))$.

Recall that $I = \langle F \rangle = \langle f_1,f_2,\ldots,f_N \rangle$,
let $J = \langle G^{(d)} \rangle \supset I$ be the ideal generated
by the polynomials in the augmented system. The reversed
containment holds for the dual spaces: $\D0[I] \supset \D0[J]$.

There is a 1-to-1 correspondence between linear differential
operators and linear differential functionals:
\begin{equation}
\sum \lambda_\beta\bfp^\beta \longleftrightarrow \sum
\lambda_\beta\beta!\Delta_\beta.
\end{equation}
Let $\phi: \bC[\bfp] \to \D0$ and $\tau: \D0 \to \bC[\bfp]$ be the
corresponding bijections.

As in Section \ref{sec-multiplicity} we order terms $\Delta_\beta$
with $\succeq$, a global monomial order. Notice that since the
choice of coefficients of the operator $Q$ is generic, $\beta_Q =
\IN_\succeq(Q) = \IN_\succeq(\phi(Q))$ is the maximal element of
the set $\IN_\succeq(\D0[I])\cap\{\beta : |\beta|\leq d \}$.

Next, we use the condition that $f_i$ form a standard basis. Since
the corners of the staircase correspond to the initial terms of
$f_i$, by Lemma \ref{lem-drop} the staircase created with the
corners at $\IN_\geq(Q\cdot(\x^\alpha f_i))$ bounds the set
$\{\beta-\beta_Q\ |\ \beta\in\IN_\succeq(\D0[I])\} \cap \ZZ_{\geq
0}^n$, which, therefore, contains the initial support of $\D0[J]$.~\qed

\begin{corollary}
If there exist a local monomial order $\geq$ such that the minimal
(standard) monomial in the set $\{\x^\alpha \notin \IN_\geq(I) :
|\alpha|\leq d \}$ is minimal in the set of all standard
monomials, then $\x^*$ is a regular solution of $G^{(d)}(\x) =
\zero$.
\end{corollary}

Ideally we would like to be able to drop the assumption of the
original polynomials forming a standard basis, since computing
such a basis is a complex symbolic task, whereas our interest lies
in the further numericalization of the approach. The following
weaker statement works around this restriction.

Assuming $\x^*=\zero$, let $\supp(F) = \bigcup_{j=1,2,\ldots,N}
\supp(f_j)$.

\begin{proposition}\label{prop-corank-drop}
Assume $A(\zero) = \zero$. Let  $d_0 = \min\{|\alpha| : x^\alpha
\in \supp(F) \}$.

Then, in the notation of Theorem \ref{thm-fixed-lambda-general},
for a generic deflating operator $Q$ the system $G^{(d)}(\x) =
\zero$, where, $d<d_0$ is a deflation of the original system
$F(\x)=\zero$ at the origin.

Moreover, if $d = d_0-1$ then the Jacobian of $G^{(d)}(\zero)$ is
not equal to zero.
\end{proposition}

\noindent {\em Proof.}
Fix a local monomial ordering that respects the degree. With the
above assumptions, the initial ideal $\IN (\langle F \rangle)$
will contain monomials of degree at least $d_0$. On the other
hand, for a generic choice of the deflating operator $Q$ the
support $\supp(G^{(d)})$ would contain a monomial of degree less
than $d_0$. Therefore, there exists a monomial in
$\IN(\langle\supp(G^{(d)})\rangle)$ that is not in $\IN (\langle F
\rangle)$, hence, $G^{(d)}$ is a deflation.

If $d = |d_0|-1$, then there is such monomial of degree 1, which
means that the Jacobian of the augmented system is nonzero.~\qed

\begin{remark}{\rm \label{rem-homogeneous}
  Note that if the deflation order $d$ is as in
  Proposition \ref{prop-corank-drop}, then it suffices to take an
  arbitrary \emph{homogeneous} deflation operator of order $d$.
}
\end{remark}

Next we explain the practical value of
Proposition~\ref{prop-corank-drop}. Let $K = \ker A(\zero)$ and $c
= \corank A(\zero) = \dim K$. Without a loss of generality we may
assume that $K$ is the subspace of $\bC^n$ has
$\{x_1,\ldots,x_c\}$ as coordinates.

Now consider the system $F'(x_1,\ldots,x_c) =
F(x_1,\ldots,x_c,0,\ldots,0)$. This system has an isolated
solution at the origin, and Proposition \ref{prop-corank-drop} is
applicable, since the Jacobian is zero. Moreover, if we take the
deflation of order $d=d_0-1$ of the original system $F$, with
$d_0$ coming from the Proposition, the corank of the Jacobian the
augmented system $G^{(d)}$ is guaranteed to be lower than that
of~$A(\zero)$.

\medskip

Let us go back to the general setup: an arbitrary isolated
solution $\x^*$, the Jacobian $A(\x^*)$ with a proper kernel $K$,
etc. Algorithm~\ref{alg_MinOrderForCorankDrop}
is a practical algorithm that can be executed numerically knowing
only an approximation to~$\x^*$.

\noindent {\em Proof.} [Proof of correctness of
Algorithm~\ref{alg_MinOrderForCorankDrop} for $\x^0 = \x^*$.]
We can get to the special setting
of Proposition \ref{prop-corank-drop} in two steps. First, apply
an affine transformation that takes $\x^*$ to the origin and $\ker
A(\x^*)$ to the subspace $K$ of $\bC^n$ spanned by the first
$c=\corank A(\x^*)$ standard basis vectors. Second, make a new
system $F'(x_1,\ldots,x_c)=0$ by substituting the $x_i=0$ in $F$
for $i>c$.

Let $\gamma' =(\gamma'_1,\ldots,\gamma'_c)\in K$ be the image of
the generic vector $\gamma$ under the linear part of the affine
transform. Then $H(t)=F'(\gamma'_1 t,\ldots,\gamma'_c t)$.

Since $\gamma'$ is generic, the lowest degree $d_0$ of the
monomial in $\supp(F')$ is equal to $\min\{ a\ |\ t^a \in \supp
H(t) \}$. According to the Proposition \ref{prop-corank-drop} and
the discussion that followed, $d=d_0-1$ is the minimal order of
deflation that will reduce the rank of the system.~\qed

\begin{remark}{\rm \label{rem-truncated-deflation-matrix}
  In view of Remark \ref{rem-homogeneous} it would be enough to use any
  \emph{homogeneous} deflation operator of order $d$:
\begin{equation}
    Q = \sum_{|\beta| = d}
        \lambda_\beta\partial^\beta \ \in \bC[\bfp],
\end{equation}
such that the vector $\bfl$ of its coefficients is in the kernel of the
\emph{truncated deflation matrix}, which contains
only the rows corresponding to the original polynomials $F$ and only
the columns labelled with $\partial^\beta$ with $|\beta|=d$.
}
\end{remark}

\subsection{Indeterminate Multipliers}

As in Section \ref{sec-old-deflation}, we now consider
indeterminate $\lambda_\beta$. Now we should think of the
differential operator $L(\bfl) \in \bC[\bfl, \partial]$ and of
additional equations $g_{j,\alpha}(\x, \bfl) \in \bC[\x,\bfl]$ as
depending on $\bfl$.

\noindent {\em Proof of Theorem~\ref{thm-deflation-step-general}.}
Picking $m = \corank(A^{(d)}(\x^*))$ generic linear equations
$h_k$ guarantees that for $\x = \x^*$ the solution for $\bfl$
exists and is unique; therefore, the first part of the statement
is proved.

The argument for the drop in the multiplicity is similar to that
of the proof of Theorem \ref{thm-deflation-step}.~\qed

\section{Computational Experiments}

We have implemented our new deflation methods in
PHCpack~\cite{V99} and Maple.  Below we report on two examples.
 
One crucial decision in the deflation algorithm is the determination
of the numerical rank, for which we may use SVD or QR in the rank-revealing
algorithms.  Both SVD and QR are numerically stable.  We summarize the
result from~\cite[page~118]{Dem97} for the problem of solving an
overdetermined linear system $A \x = b$.  The solution obtained by
QR or SVD minimizes the residual
$|| (A + \delta A) \widetilde{\x} - (b + \delta b) ||_2$ where
the relative errors have the same magnitude as the machine
precision~$\epsilon$:
\begin{equation}
   {\rm max}\left( \frac{||\delta A||_2}{||A||_2} ,
                   \frac{||\delta b||_2}{||b||_2} \right) = O(\epsilon).
\end{equation}

\subsection{Running Example 1}

To find initial approximations for the roots
of the system~(\ref{equ-running}), we must first make the system
``square'', i.e.: having as many equations as unknowns, so
we may apply the homotopies available in PHCpack~\cite{V99}.
Using the embedding technique of~\cite{SV00}
(see also~\cite{SW05}), we add one
slack variable~$z$ to each equation of the system,
multiplied by random complex constants
$\gamma_1$, $\gamma_2$, and $\gamma_3$:

\begin{equation} \label{equ-embed-running}
  E(\x,z)
  = \left\{
       \begin{array}{c}
          x_1^3 + x_1 x_2^2 + \gamma_1 z = 0 \\
          x_1 x_2^2 + x_2^3 + \gamma_2 z = 0 \\
          x_1^2x_2 + x_1x_2^2 + \gamma_3 z = 0.
       \end{array}
    \right.
\end{equation}
Observe that the solutions of the original system $F(\x) = \zero$
occur as solutions of the embedded system $E(\x,z) = \zero$
with slack variable~$z = 0$.
At the end points of the solution paths defined by a homotopy
to solve $E(\x,z) = \zero$, we find nine zeroes close to the origin.
These nine approximate zeroes are the input to our deflation algorithm.

The application of our first deflation algorithm in~\cite{LVZ06}
requires two stages.
The Jacobian matrix of $F(\x) = \zero$ has rank zero at $(0,0)$.
After the first deflation with one multiplier, the rank of the
Jacobian matrix of the augmented system $G(\x,\bfl_1) = \zero$
equals one, so the second deflation step uses two multipliers.
After the second deflation step, the Jacobian matrix has full
rank, and~$(0,0)$ has then become a regular solution.
Newton's method on the final system then converges again quadratically
and the solution can be approximated efficiently with great accuracy.
Once the precise location of a multiple root is known,
we are interested in its multiplicity.
The algorithm of~\cite{DZ05} reveals that the multiplicity
of the isolated root equals seven.

Starting at a root of low accuracy, at a distance of $10^{-5}$
from the exact root, the numerical implementation of
Algorithm~\ref{alg_MinOrderForCorankDrop} predicts two as the
order, using $10^{-4}$ as the tolerance for the vanishing of the
coefficients in the univariate interpolating polynomial.  The
Jacobian matrix of the augmented system $G^{(2)}$ has full rank so
that a couple of iterations suffice to compute the root very
accurately.

\subsection{A Larger Example}

The following system is copied from~\cite{Lec02}:
\begin{equation}
   F(\x) =
   \left\{
     \begin{array}{c}
        2 x_1 + 2 x_1^2 + 2 x_2 + 2 x_2^2 + x_3^2 - 1 = 0 \\
        (x_1 + x_2 - x_3 - 1)^3 - x_1^3 = 0 \\
        (2 x_1^3 + 2 x_2^2 + 10 x_3 + 5 x_3^2 +5)^3-1000 x_1^5 = 0.
     \end{array}
   \right.
\end{equation}
Counted with multiplicities, the system has 54 isolated solutions.
We focus on the solution $(0,0,-1)$ which occurs with multiplicity~18.

Although Algorithm~1 
suggests that the first-order deflation would already lower the
corank of the system, we would like to search for a homogeneous
deflation operator $Q$ of order two.

To this end we construct the (truncated) deflation matrix $\bar
A(x_1,x_2,x_3)$ with has 12 rows and only 6 columns, which
correspond to $\{\partial_1^2, \partial_1\partial_2,
\partial_1\partial_3,
\partial_2^2, \partial_2\partial_3, \partial_3^2\}$.

The kernel of $\bar A(0,0,-1)$ is spanned 
by $(1, 6, 8, -3, 0, 4)^T$ and $(0, 3, 3, -1, 1, 2)^T$.
The operator corresponding to the former,
\begin{equation}
    Q = \partial_1^2 + 6 \partial_1\partial_2
     + 8 \partial_1\partial_3 - 3 \partial_2^2 + 4 \partial_3^2,
\end{equation}
regularizes the system, since the equations
\begin{equation}
  \left\{
     \begin{array}{rcrcrcrcrcl}
       Q\cdot(x_1f_1) &=& 8 x_1 & + & 24 x_2 & + & 16 x_3 & + & 16 & = & 0 \\
       Q\cdot(x_2f_1) &=& 24 x_1 & - & 24 x_2 &   &        &   &    & = & 0 \\
       Q\cdot(x_3f_1) &=& 32 x_1 &   &        & + & 16 x_3 & + & 16 & = & 0.
     \end{array}
  \right.
\end{equation}
augmented to the original equations, give a system with the
full-rank Jacobian matrix at $(0,0,-1)$.

\section{Conclusion}

In this paper we have described two methods of computing the
multiplicity structure at isolated solutions of polynomial
systems.  We have developed a higher-order deflation algorithm that
reduces the multiplicity faster than the first-order deflation 
in~\cite{LVZ06}.

In our opinion, one of the main benefits of the higher order deflation
for the numerical algebraic geometry algorithms is the possibility to 
regularize the system in a single step. 
For that one has to determine the minimal order of such a deflation or,
even better, construct a sparse ansatz for its deflation operator.
Predicting these numerically could be a very challenging task, 
which should be explored in the future.
   
\bibliographystyle{plain}
\bibliography{hod}

\end{document}